\documentclass[10.5pt, a4paper,  reqno, english]{amsart}%draft%{article}
\usepackage[all]{xy}
\usepackage[utf8]{inputenc}
\usepackage{amssymb,amsmath,amsthm}
\numberwithin{equation}{section}
\usepackage{amsfonts}
\usepackage{bbm}
\usepackage{mathrsfs}
\usepackage{lmodern}
\usepackage{enumitem}
\usepackage{xcolor}
\usepackage{comment}

\usepackage[margin=1.5in]{geometry}

\newtheorem{thm}{Theorem}
\newtheorem{prop}{Proposition}

\begin{document}
\title{A note on Benli and Pollack's paper on small prime power residues}
\author{Crystel Bujold}
\address{Département de mathématiques et de statistique, Université de Montréal, CP 6128
succ. Centre-Ville, Montréal, QC H3C 3J7, Canada}
\email{bujoldc@dms.umontreal.ca}
\begin{abstract}
    In this note, we revisit a result of Benli's and Pollack's on the number of small prime $k^{th}$ power residues. The proof is based on their idea of using reciprocity laws, but the argument is simplified and we prove a slightly stronger bound.
\end{abstract}

\maketitle

\section{Introduction}

The problem of finding a bound for the least $k^{th}$ power residue modulo $q$ has been a question of interest for many years. In general, finding this bound is a difficult problem that relies on a careful analysis of character sums and the associated L-functions and so, although we expect the least residues to be of size $O(q^{\epsilon})$, the best unconditional bound we have at hand for the least $k^{th}$ power residue is $O(q^{\frac{k-1}{4}+\epsilon})$ (see \cite{ell}, \cite{vinlin}, \cite{soren} ). Recently, Pollack showed \cite{PollackPNR} that in the case of quadratic residues, there are in fact many prime residues below $q^{1/4-\epsilon}$ when the modulus $q$ is a prime. This opens the interesting question of investigating the number of prime residues (and non-residues) up to a given point.  

In a recent paper, Benli and Pollack observed that, among all the results on residues, it hasn't been widely noticed that Nagell (1952, \cite{Nagell}) improved on the known bounds for the least cubic residue and that the techniques used in the proof, namely with the reciprocity laws, deserved to be further explored. Hence, in \cite{benliSPR}, exploiting Nagell's ideas, they provide lower bounds for the number of quadratic, cubic and biquadratic residues up $q^{1/2+\epsilon}$. In a subsequent paper \cite{MR4127826}, Benli improved the result, giving a lower bound for the number of $k^{th}$ power prime residues modulo a prime $p$, for all $k\geq 2$, up to $q^{\frac{k+1}{4}+\epsilon}$. 

The focus of the present note is to give a simplified proof and slightly stronger version of the results in \cite{benliSPR}. Besides using reciprocity laws, Benli's and Pollack's proof rely on sieve methods which are efficient but hardly enlightening and it is possible to avoid  them by using elementary summation techniques. In the following sections, we prove the following theorems:

\begin{thm}
Let $\epsilon \leq \frac{1}{5}$ and let $p\geq e^{312}$, then

\begin{align*}
		\#  \left\{p^{\frac{\epsilon}{2}} < q \leq 2p^{\frac{1}{2}+\epsilon} :\textrm{q is a quadratic residue}\pmod p \right\}
        & \geq \frac{1}{4} p^{\frac{\epsilon}{2}}.
\end{align*}

\end{thm}
We then prove a similar theorem for cubic residues.

\begin{thm}
Let $\epsilon \leq \frac{1}{2}$ and let $p\geq \max\{75^{\frac{\epsilon}{3}},e^{100}\}$. Then 

\begin{equation*}
\#\left\{p^{\frac{\epsilon}{3}} < q \leq 72p^{\frac{1}{2}+\epsilon} : \textrm{q is a cubic residue} \pmod p\right\}\geq \frac{1}{15}p^{\frac{\epsilon}{3}}.
\end{equation*}

\end{thm}
Finally, we show 

\begin{thm}
Let $\epsilon \leq \frac{1}{2}$ and let $p\geq \max\{60^{\frac{\epsilon}{4}},e^{85}\}$. Then

\begin{equation*}
\# \left\{x\leq q\leq 42p^{\frac{1}{2}+\epsilon}: \textrm{q is a biquadratic residue} \pmod p\right\}\geq \frac{1}{16}p^{\frac{\epsilon}{4}}.
\end{equation*}

\end{thm}

\section{Quadratic residues}

The strategy follows Pollack's ideas. Given a prime $p$, we wish to find a lower bound for the number of prime residues modulo $p$. To do so, we recall that from quadratic reciprocity
\begin{equation*}
	\left(\frac{q}{p}\right) = \left\{ \begin{array}{lr}
									\left(\frac{p}{q}\right) & \textrm{ if } p\equiv 1 \mod 4\\
                                    \left(\frac{-p}{q}\right) & \textrm{ if } p\equiv 3 \mod 4 \\                                           \end{array}    \right.                               
\end{equation*}
That is, using quadratic polynomials, we want to show that there are many primes $q$ for which $p$ or $-p$ is a square $\mod q$, according whether $p \equiv 1 \textrm{ or } 3 \mod 4$.

\subsection{The case $p \equiv 1 \mod 4$}

Let $p$ be a prime such that $p\equiv 1 \mod 4$. Following Pollack's idea, define $r=\lfloor\sqrt{p}\rfloor $ and let 
\begin{equation*}
	f(n) = (n+r)^2 -p.
\end{equation*}
Observe that $f(n)$ doesn't have fixed divisors and that it has at most 2 roots modulo every prime $q$. By applying Hansel's lemma, we also get that $f(n)$ has at most 2 roots modulo every prime power. Morevover, q divides $f(n)$ if and only if $p\equiv (n+r)^2 \mod q$. Hence, we want to find a lower bound for the number of primes $q$ dividing $f(n)$ for some integer $n$.
%%%Notice that $q^k|f(n)$ if and only if $(n+r)^2 -p\equiv 0 \mod q^k$, and since $q$ doesn't divide $f'(n)$, by Hensel's lemma this also has at most 2 solutions modulo $q^k$ and therefore, there are at most $\frac{2x}{q^k}$ possible values of $n$ for which this holds if $q \leq x^{\frac{1}{k}}$ and only 2 possible values of $n$ if $q > x^{\frac{1}{k}}$. %%%

Let $x= p^{\epsilon}$, and consider

\begin{align*}
	\sum_{n\leq x} \log(f(n)) &= \sum_{n\leq x}\sum_{k\geq 1} \sum_{q^k|f(n)} \log q\\
    						  &= \sum_{q\leq x} \sum_{k\geq 1}\sum_{\substack{n\leq x \\ q^k|f(n)}}\log q +    						                      \sum_{x< q\leq2r x} \sum_{k\leq \frac{\log (2rx)}{\log q}}\sum_{\substack{n\leq x \\ q^k|f(n)}}\log q. \\
\end{align*}
For the first sum, we have

\begin{align*}
	\sum_{q\leq x} \sum_{k \geq} \sum_{\substack{n\leq x \\ q^k | f(n)}} \log q &\leq 2x\sum_{q\leq x}\sum_{k\geq 1} \frac{\log q}{q^k}\\
     &= 2x \sum_{q\leq x} \log q \sum_{k \geq 1} \left(\frac{1}{q}\right)^k\\
     &= 2x \sum_{q \leq x} \frac{\log q}{q}\left(\frac{1}{1- \frac{1}{q}}\right)\\
     &=2x\sum_{q\leq x} \frac{\log q}{q-1}\\
     & \leq 2xc_1(2 + \log x)\\
     & = 2c_1\left( \frac{2}{\log p} + \epsilon \right)x\log p\\
\end{align*}
where $c_1 = \left(1+ \frac{0.15}{\log^3 x}\right)$.\\
Now from the second sum, we get

\begin{align*}
	\sum_{x < q \leq 2rx} \sum_{k \leq \frac{\log (2rx)}{\log q}} \sum_{\substack{n\leq x \\ q^k | f(n)}} \log q 
    &\leq \sum_{x < q \leq 2rx} \sum_{k \leq \frac{\log (2rx)}{\log q}} \sum_{\substack{n\leq x \\ q | f(n)}} \log q \\
    & = \sum_{x < q\leq 2rx}\sum_{\substack{n \leq x \\  q|f(n)}} \frac{\log (2rx)}{\log q}\log q\\
    & \leq 2\log(2rx) \cdot \#\{x < q \leq 2rx : q| f(n), \textrm{ for some } n\}\\
    & = \left(1+2\epsilon +\frac{\log 4}{\log p}\right) \cdot \#\{x < q \leq 2rx : q| f(n), \textrm{ for some } n\}\\
\end{align*}
On the other hand, using stirling's approximation, given $\epsilon \geq \frac{\log\log p}{\log p}$ we have

\begin{align*}
	\sum_{n\leq x} \log(f(n)) &= \sum_{n\leq x} \log((n+r)^2-p)\\
    						  &= \sum_{n\leq x} \log (2nr) + \log \left( 1+ \frac{n^2 + r^2 -p}{2nr}\right) \\
                              &\geq \sum_{n\leq x}\log 2 + \log n + \log r + \log\left(1-\frac{1}{n}\right)\\
                              &\geq (x-1)\log 2 +\log (\lfloor x \rfloor !) +  (x-1)\log r -(x+2\log x)\\
                              &\geq  x\log x + x\log r -3x\\
                              &= \left(\frac{1}{2} + \epsilon - \frac{3}{\log p}\right)x\log p
\end{align*}
Hence, altogether we get

\begin{equation*}
	\left(\frac{1}{2} + \epsilon - \frac{3}{\log p}\right)x\log p \leq    2c_1\left( \frac{2}{\log p} + \epsilon \right)x\log p + \left(1+2\epsilon +\frac{\log 4}{\log p}\right)\log p \cdot \#\{x < q \leq 2rx : q| f(n), \textrm{ for some } n\leq x \}
\end{equation*}
so that

\begin{equation*}
	 \#\{x < q \leq 2rx : q| f(n), \textrm{ for some } n\leq x\}\geq 	\frac{\left(\frac{1}{2} + (1-2c_1)\epsilon - \frac{4(1+c_1)}{\log p}\right)x}{\left(1+2\epsilon +\frac{\log 4}{\log p}\right)}
\end{equation*}
Therefore we deduce that, 

\begin{equation*}
	 \#  \left\{p^{\epsilon} < q \leq 2p^{\frac{1}{2}+\epsilon} : q| f(n) \textrm{ for some } n\leq p^{\epsilon} \right\} \geq   \frac{1}{4}p^{\epsilon},
\end{equation*}
given that $\epsilon \leq \frac{1}{4}$ and $ p \geq e^{144}$. 

\subsection{The case $p\equiv 3 \mod 4$}

Let $p$ be a prime such that $p \equiv 3 \mod 4$. Then by the law of quadratic reciprocity, we know that for a prime $q$,  $\left(\frac{q}{p}\right) = \left(\frac{-p}{q}\right)$, i.e. $q$ is a quadratic residue modulo $p$ if and only if $-p$ is a square modulo $q$.

Now recall the following proposition,

\begin{prop}	An integer $m$ is properly represented by a binary quadratic form of discriminant $D$ if and only if $D$ is a square modulo $4m$.
\end{prop}

Hence, letting $D = -p$, for any integer $m$ properly represented by a binary quadratic form of discriminant $D$, it follows that for, all prime $q|m$, then $D$ is a square modulo $q$ . By quadratic reciprocity, we will get that these primes $q$ are quadratic residues modulo $p$.

We follow essentially the same lines as before except that we have to work a little harder to find the appropriate quadratic form giving us the lower bound we are looking for. We take $f(n)= an^2+bn+c$ to be defined as in Pollack's paper, so that $p^{\frac{1}{2}-\frac{\epsilon}{2}} < a < p^{\frac{1}{2}}$, $|b|\leq a$ and $p^{\frac{1}{2}}\leq c \leq p^{\frac{1}{2}+ \frac{\epsilon}{2}}$. Observe that in that case, taking $x = p^{\frac{\epsilon}{2}}$, if $q|f(n)$ for some $n\leq x$ then $q \leq 2p^{\frac{1}{2}+\epsilon}$.

From there, the argument is essentially identical as for the case $p\equiv 1 \mod 4$. We have on one hand, replacing $2rx$ by its value $2p^{\frac{1}{2} + \epsilon}$,

\begin{align*}
	\sum_{n\leq x} \log(f(n)) &= \sum_{n\leq x}\sum_{k\geq 1} \sum_{q^k|f(n)} \log q\\
        					 &= \sum_{q\leq x} \sum_{k\geq 1}\sum_{\substack{n\leq x \\ q^k|f(n)}}\log q +    						                      \sum_{x< q\leq 2p^{\frac{1}{2} + \epsilon}} \sum_{k\leq \frac{\log (2p^{\frac{1}{2} + \epsilon})}{\log q}}\sum_{\substack{n\leq x \\ q^k|f(n)}}\log q. \\
                             & \leq 2c_1\left( \frac{2}{\log p} +\frac{\epsilon}{2} \right)x\log p + \left(1+\epsilon +\frac{\log 4}{\log p}\right) \log p \cdot \#\{x < q \leq 2p^{\frac{1}{2}+\epsilon} : q| f(n), \textrm{ for some } n\}\\
\end{align*}
On the other hand, we have 

\begin{align*}
	\sum_{n\leq x} \log(f(n)) &= \sum_{n \leq x} \log(an^2 +bn +c)\\
   			 &\geq \sum_{n \leq x} \log(an^2 +bn )\\
    		& = \sum _{n \leq x} \log(a) + 2\log n + \log \left(1-\frac{1}{n}\right)\\
            & \geq (x-1)\log a +2\log (\lfloor x \rfloor!) - \sum_{n\leq x}\frac{1}{n-1} \\
            & \geq x\log a +2(x\log x-x-2\log x) -\log x + \log a\\
            & \geq \left(\frac{1}{2} + \frac{\epsilon}{2}\right)x\log p -4x\\
            &\geq \left( \frac{1}{2} + \frac{\epsilon}{2} - \frac{4}{\log p}\right) x\log p\\
\end{align*}
Hence we get

\begin{equation*}
	%\begin{split}
		c_1\left(\frac{4}{\log p} + \epsilon\right)x\log p +  \left(1 + \epsilon +\frac{\log 4}{\log p} \right) \log p \cdot \# \small{ \left\{x < q \leq 2p^{\frac{1}{2}+ \epsilon} : q| f(n) \textrm{ for some } n\leq x \right\} }
        \geq \left( \frac{1}{2} - \frac{\epsilon}{2} - \frac{4}{\log p} \right) x \log p 
   % \end{split}
\end{equation*}
and therefore

\begin{align*}
		\#  \left\{p^{\frac{\epsilon}{2}} < q \leq 2p^{\frac{1}{2}+\epsilon} : q| f(n) \textrm{ for some } n\leq p^{\frac{\epsilon}{2}} \right\} &\geq \frac{\left(\frac{1}{2}-\frac{5\epsilon}{6} -\frac{9}{\log p}\right)}{\left(1 + \epsilon +\frac{\log 4}{\log p}\right)}x \\
        & \geq \frac{1}{4} p^{\frac{\epsilon}{2}}
\end{align*}
for $\epsilon \leq \frac{1}{5}$ and $p\geq e^{312}$.\\

Observe that if we take the argument with $f(n) = x^2 + x + \frac{1+p}{4}$ and $q\leq p$, i.e. the special case $\epsilon = \frac{1}{2}$, then taking $x=p^{\delta}$ for $\delta < \frac{1}{2}$ we get

\begin{equation*}
		\#  \left\{p^{\delta} < q \leq p : q| f(n) \textrm{ for some } n\leq p \right\} \geq (1-2\delta) p^{\delta} + O(x).
\end{equation*}
%%% We can get a small gain on the range for epsilon and p by summing over y and y in the quadratic form%%%

\section{Cubic residues}

The argument for the cubic case is essentially the same as for the quadratic case, except that we'll get our cubic polynomial from the following proposition.

\begin{prop}
	Let $p\equiv 1 \mod 3 $, let $L$ and $M$ be the unique integers such thath $4p = L^2 + 27M^2$. Say $q \neq 2,3,p$, then $q$ is a cubic residue $\mod p$ if and only if
     \begin{equation*}
     	\frac{L}{3M}\equiv \frac{x^3-9x}{3(x^2-1)} \mod q,
     \end{equation*}
     for some $x\in \mathbb{Z}$.
\end{prop}

Now if $3L(x^2-1)-3M(x^3-9x) \equiv 0 \mod q$, then the above condition is satisfied and $q$ is a cubic residue so that we define 

\begin{equation*}
f_0(n) = L(n^2-1)-M(n^3-9n) 
\end{equation*}
and in order to have the polynomial being greater than 1 for all values of $n$, we will use
\begin{equation*}
	f(n) = f_0(3n)=L(9n^2-1)-27M(n^3-n).
\end{equation*}

For the cubic case, we let $\epsilon \leq \frac{1}{2}$ and $x= p^{\epsilon/3}$,  so that if $q|f(n)$, then $q \leq 72p^{\frac{1}{2} + \epsilon}$.The first part of the argument is identical to the first part in the quadratic case, except that now the polynomial can have at most 3 roots modulo every prime power. Hence we have that

\begin{align*}
		\sum_{n\leq x} \log(f(n)) &= \sum_{n\leq x}\sum_{k\geq 1} \sum_{q^k|f(n)} \log q\\
        					 &= \sum_{q\leq x} \sum_{k\geq 1}\sum_{\substack{n\leq x \\ q^k|f(n)}}\log q + \sum_{x< q\leq 2p^{\frac{1}{2} + \epsilon}} \sum_{k\leq \frac{\log (72p^{\frac{1}{2} + \epsilon})}{\log q}}\sum_{\substack{n\leq x \\ q^k|f(n)}}\log q. \\
                             & \leq 3c_1\left( \frac{2}{\log p} +\frac{\epsilon}{3} \right)x\log p + 3\left(\frac{1}{2}+\epsilon +\frac{\log 72}{\log p}\right) \log p \cdot \#\{x < q \leq 72p^{\frac{1}{2}+\epsilon} : q| f(n), \textrm{ for some } n\}\\
\end{align*}
On the other hand, observe that $\max\{L,M\}\geq \sqrt{\frac{p}{7}}$, so we first suppose that $\max\{L,M\} = M$. That is, we have 

\begin{align*}
\sum_{n\leq x} \log(f(n)) &\geq \sum_{n\leq x} \log (8+ 27M(n^3-n))\\
&\geq \log 8 + \sum_{2\leq n\leq x} \log(27M(n^3-n))\\
&\geq \log 8 + (x-2)\log(27M) + 3\log (\lfloor x \rfloor!) -\log 8 -2\sum_{2\leq n\leq x}\log\left(1-\frac{1}{n^2}\right)\\
& \geq (x-2)\log\left(27\left(\frac{p}{7}\right)^{1/2}\right) + 3(x\log x -x-2\log x) - 2\sum_{n\leq x} \frac{1}{n^2}\\
& \geq \left(\frac{1}{2} + \epsilon - \frac{0.68}{\log p} - \frac{3}{x}\right) x\log p.\\
\end{align*}

Therefore, putting it all together as in the quadratic case, we obtain, say for $x\geq 50$ and $p \geq e^{60}$,

\begin{equation*}
\#\{x < q \leq 72p^{\frac{1}{2}+\epsilon} : q| f(n), \textrm{ for some } n\}\geq \frac{1}{10}p^{\epsilon/3}.
\end{equation*}

Now if we repeat the argument with $\max{M, L} = L$, we get
\begin{align*}
\sum_{n\leq x} \log(f(n)) &\geq \sum_{n\leq x} \log(L(9n^2-1))\\
& \geq \left(\frac{1}{2} + \frac{2\epsilon}{3} - \frac{0.78}{\log p} - \frac{4}{x}\right) x\log p.
\end{align*}
From which we deduce that, if $x\geq 75$ and $p \geq e^{100}$, then
\begin{equation*}
\#\{x < q \leq 72p^{\frac{1}{2}+\epsilon} : q| f(n), \textrm{ for some } n\}\geq \frac{1}{15}p^{\epsilon/3}.
\end{equation*}

\section{Biquadratic residues}

Just as in the previous cases, we use an explicit criterion for biquadratic reciprocity.

\begin{prop}
Let $p\equiv 1 \mod 4$ and write $p$ with its unique representation $p = L^2 + 4M^2$, and let $q* = (-1)^{(q-1)/2}q$, then $q*$ is a biquadratic residue if and only if 
\begin{equation*}
\frac{L}{2M} \equiv \frac{x^4-6^2+1}{4(x^3-x)} \mod q
\end{equation*}
for some $x\in \mathbb{Z}$.
\end{prop}

In an analogous way to the previous cases, we set $f_0(x) = 2L(27x^3-3x) + M(x^4-6x+1)$ and we define $f(n) = f_0 (3n) = 6L(9n^3-n) + M(81n^4-54n^2 +1) > 0$ for all $n\geq 1$. Once again, we observe that if a prime $q|f(n)$ then the condition for reciprocity is satisfied and $q$ is a biquadratic residue modulo $p$. Therefore, we wish to find a lower bound on the number of primes dividing $f(n)$.

The first part of the argument is again identical to the quadratic case, with a maximum of 4 roots modulo every prime power instead of 2. That is, taking $x=p^{\epsilon/4}$ so that for $n\leq x$, $f(n) \leq 42p^{\frac{1}{2} + \epsilon}$, we have

\begin{equation*}
\sum_{n\leq x} \log (f(n)) \leq 4c_1 \left(\frac{2}{\log p}+ \frac{\epsilon}{4}\right) x\log p + 4\log(42p^{\frac{1}{2}+\epsilon}) \cdot \# \{x\leq q\leq 42p^{\frac{1}{2}+\epsilon}: q|f(n) \textrm{ some } n\geq 1\}.
\end{equation*}
On the other hand we have 

\begin{equation*}
\sum_{n\leq x} \log (f(n)) = \sum_{n\leq x} \log (6L(9n^3-n) + M(81n^4 -54n^2 +1)).
\end{equation*}
Observe that from the equality $p = L^2 + 4M^2$, we have tha $\max\{L,M\}\geq \sqrt{\frac{p}{5}}$. Morover, notice that since $6(9n^3-n) \leq (81n^4-54^2+1)$ for all $n\leq 1$, it follows that $f(n) \geq 6\sqrt{\frac{p}{5}}(9n^3-n)$. Hence

\begin{align*}
\sum_{n\leq x} \log (f(n)) &\geq \sum_{n\leq x} \log(63Ln^3) + \sum_{n\leq x} \log\left(1-\frac{1}{9n^2}\right)\\
&\geq \left( \frac{1}{2} + \frac{3\epsilon}{4} + \frac{1}{3\log p} - \frac{5}{4x}\right) x\log p.
\end{align*}
Therefore, putting both bounds together, we deduce as before that for $x\geq 60$ and $p\geq e^{84}$,

\begin{equation*}
\# \{x\leq q\leq 42p^{\frac{1}{2}+\epsilon}: q|f(n) \textrm{ some } n\geq 1\}\geq \frac{1}{16}p^{\frac{\epsilon}{4}}.
\end{equation*}

\bibliographystyle{plain}
\bibliography{biblio.bib}
%Remark: The biquadratic case only covers $p\equiv 1 \mod 4$. What about $p\equiv 3 \mod 4$?

\end{document}